\documentclass[12pt]{amsart}
\usepackage{amsfonts, amsmath, latexsym, epsfig}
\usepackage{amssymb}
\usepackage{epsf}
\usepackage{url}

\newcommand{\RR}{\ensuremath{\mathbb{R}}}

\newcommand{\ZZ}{\ensuremath{\mathbb{Z}}}
\newcommand{\TT}{\ensuremath{\mathbb{T}}}
\newtheorem{proposition}{Proposition}
\newtheorem{theorem}{Theorem}

\newtheorem{lemma}{Lemma}

\newtheorem{definition}{Definition}
\def\QuotS#1#2{\leavevmode\kern-.0em\raise.2ex\hbox{$#1$}\kern-.1em/\kern-.1em\lower.25ex\hbox{$#2$}}

\DeclareMathOperator{\Kert}{Ker}

\begin{document}

\author{Mathieu Dutour Sikiri\'c}
\address{Mathieu Dutour Sikiri\'c, Rudjer Boskovi\'c Institute, Bijenicka 54, 10000 Zagreb, Croatia, Fax: +385-1-468-0245}
\email{mdsikir@irb.hr}
\thanks{The author has been supported by the Croatian Ministry of Science, Education and Sport under contract 098-0982705-2707. The author also thanks Y. Itoh for having invited him in Hayama where this research was initiated}

\title{Plane square tilings}

\date{}

\maketitle

\begin{abstract}
We consider here square tilings of the plane.
By extending the formalism introduced in \cite{tutte} we build a
correspondence between plane maps endowed with an harmonic vector and
square tilings satisfying a condition of regularity.
In the case of periodic plane square tiling the relevant space
of harmonic vectors is actually isomorphic to the first homology
group of a torus.
So, periodic plane square tilings are described by two parameters and the
set of parameters is split into angular sectors.

The correspondence between symmetry of the square tiling and symmetry
of the plane maps and harmonic vectors is discussed and a method for
enumerating the regular periodic plane square tilings having $r$ orbits
of squares is outlined.
\end{abstract}

\section{Introduction}

\begin{figure}
\begin{center}
\epsfig{height=48mm, file=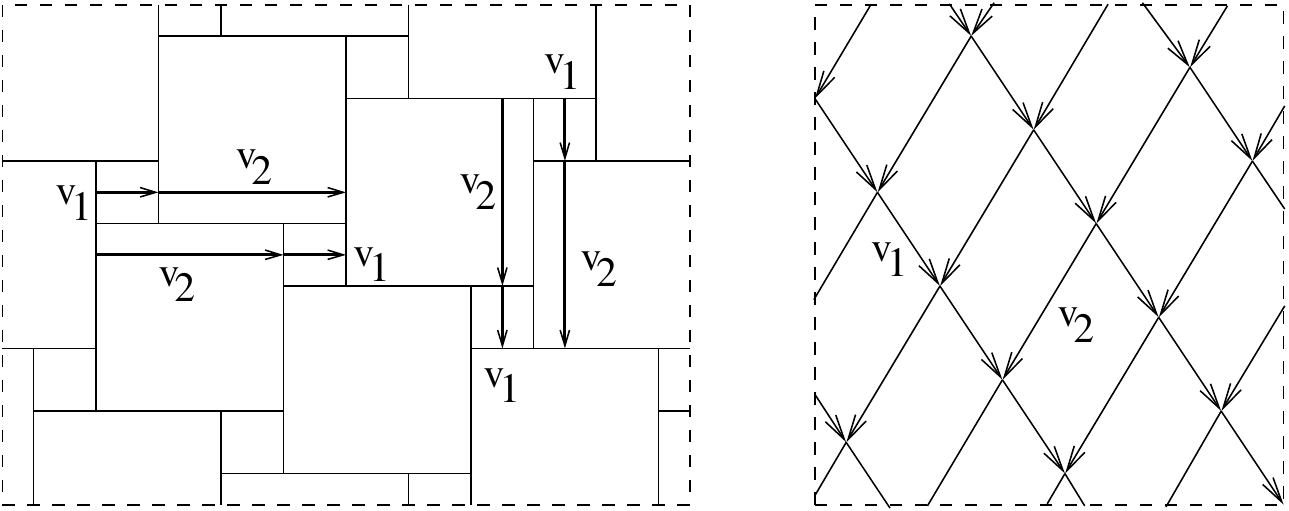}\par
\end{center}
\caption{The city square tiling and the corresponding plane map}
\label{CitySquareTiling}
\end{figure}

The face-to-face plane square tiling $[0,1]^2+\ZZ^2$ is the most basic example
of plane tiling.
An interesting variation, found on many city squares, is the plane tiling
by squares of two different size and no face-to-face adjacency
(see Figure \ref{CitySquareTiling}).
The next interesting case, i.e. the classification of plane tilings
with no face-to-face adjacency and $3$ orbits of squares of
different square sizes
was begun in \cite{grunbaum} and completed
in \cite{martini,schattschneider,bolcskei} by using classical
methods of tiling theory.
In parallel to this in \cite{tutte} an electrical network formalism
was introduced for the dissection of rectangles into squares of
unequal length
(see \cite{gambini} for a very good history of the problem).

In Section \ref{CombFormalism}, we propose a combinatorial
formalism for describing plane square tilings satisfying the condition
that there is no line containing an infinity of edges of square.
The correspondence associates to any such tiling a plane map ${\mathcal M}$
and an harmonic vector $w$ on it.
Reversely, for every plane map and harmonic vector we define an unextendible
plane square packing, which is actually a tiling in many cases.
This directly extends the work of \cite{tutte} for dissection of rectangles
to square plane tilings.

In Section \ref{HomologyFormalism}, we identify, in the periodic
case, the harmonic vectors to homology classes.
That is, to any edge finite periodic tiling ${\mathcal T}$
of the plane by squares, we associate a periodic plane maps ${\mathcal M}$.
Then the tiling ${\mathcal T}$ is described in terms of two
parameters which correspond to the first homology group of ${\mathcal M}$
quotiented by its group of translations.
However, in the non-periodic case, such an identification is not
possible and the variety of tilings is much greater.

In Section \ref{TypeSquareTiling}, we see that the space of possible
parameters of a square packing does not correspond in general to the
full homology group.
In general, the two dimensional space of the first homology group
of a torus is split into angular sector containing the origin.
We look at how the $8$ periodic square tilings of
\cite{martini,schattschneider,bolcskei} are thus described.

The wallpaper symmetry groups of such tilings are restricted and
in Section \ref{PossibleSymmetry} we consider how the symmetries of
a periodic square tiling correspond to the symmetries of the corresponding
toroidal map.
We then introduce a systematic method for enumerating the periodic regular
square tilings with $r$ orbits of squares.

\section{Harmonic vector formalism}\label{CombFormalism}

\begin{figure}
\begin{center}
\begin{minipage}[b]{3.2cm}
\centering
\epsfig{height=20mm, file=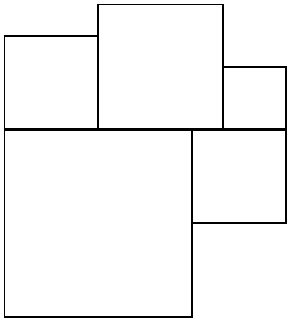}\par
\end{minipage}

\end{center}
\caption{An horizontal face}
\label{HorizFace}
\end{figure}

For a tiling of the plane or a rectangle by squares aligned along
the $x$ axis, we call {\em horizontal edge}, respectively {\em vertical
edge}, the edges of squares aligned with $x$-axis, respectively $y$-axis.
An {\em horizontal face} is the set of horizontal edges contained in
an horizontal line.
See an example on Figure \ref{HorizFace}.
Similarly we define vertical faces.

\begin{figure}
\begin{center}
\begin{minipage}[b]{5.2cm}
\centering
\epsfig{height=37mm, file=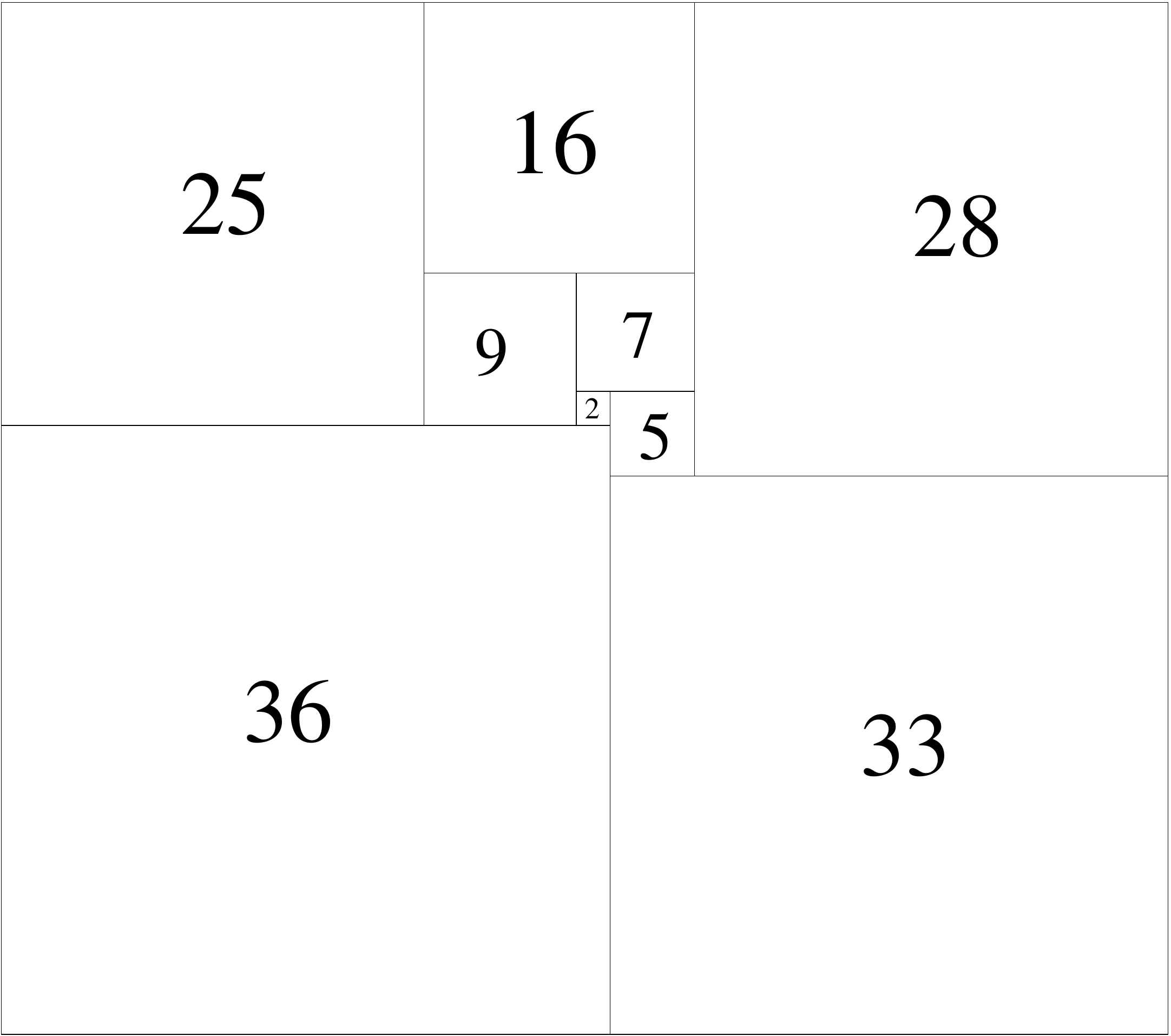}\par
\end{minipage}
\begin{minipage}[b]{5.2cm}
\centering
\epsfig{height=37mm, file=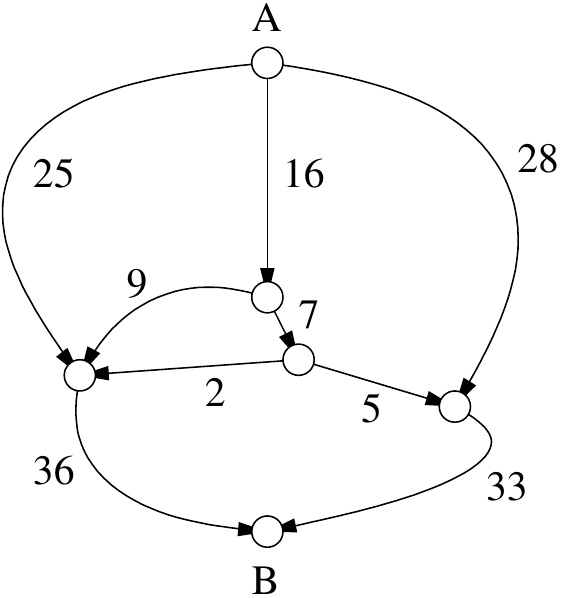}\par
\end{minipage}

\end{center}
\caption{Dissection of rectangle into squares and corresponding electrical network}
\label{DissecRectangle}
\end{figure}

In \cite{tutte} the problem of dissecting a rectangle into squares of
unequal size was considered.
There, to any such dissection into $n$ squares is associated an
electrical circuit on $n$ wires each of conductance $1$.
This electrical circuit forms a finite planar graph between between two
points $A$ and $B$. For an horizontal face, the sum of lengths of
above edges should be equal to the sum of lengths of below edges;
this is translated in conservation of currents at vertices.
Similarly equality of sum of lengths at vertical edges correspond
to voltage equality in faces.
Picture \ref{DissecRectangle} from \cite{tutte} indicates the
correspondence in an example.

For tiling of the plane by squares, the correspondence
goes along the same lines but some technical conditions are needed.
A square tiling is called {\em edge finite} if all its vertical and
horizontal faces contain a finite number of vertical and horizontal
edges.
A vertex is called {\em singular} if it is contained in exactly $4$
different squares. A plane square tiling is called {\em regular} if
it is edge finite and has no singular vertices.

By a {\em plane graph} ${\mathcal M}$ we mean a graph which can be drawn on the plane
such that any two edges if not disjoint interest only at their vertices.
The {\em faces} of ${\mathcal M}$ are the connected components of 
the plane minus vertices and edges.
We denote by $V({\mathcal M})$, $E({\mathcal M})$ and $F({\mathcal M})$ the set of vertices, edges and faces
of ${\mathcal M}$.
A {\em plane map} is a connected plane graph such that any face is bounded
and any compact subset of the plane contains only a finite number of
vertices and edges. Such a map is necessarily infinite.
A plane map is called {\em periodic} if it is invariant under translation
in two non-parallel directions.
Of course periodic plane maps correspond to toroidal maps
and we can use the notions of vertices, edges and faces for tori.

Given a plane map ${\mathcal M}$ a {\em directed edge} $\overrightarrow{e}$ is a pair $(v,e)$
with $v$ a vertex and $e$ an edge containing $v$.
For a vertex $v$, edge $e$, or face $f$ we denote by $\overrightarrow{E}(v)$, $\overrightarrow{E}(e)$, or $\overrightarrow{E}(f)$ the set
of directed edges originating from $v$, contained in $e$,
or having $f$ on the right.
An {\em harmonic vector} $w=(w_{\overrightarrow{e}})$ of a plane map ${\mathcal M}$ is
a vector defined over the set of directed edges such that for any $v\in V({\mathcal M})$,
$e\in E({\mathcal M})$, $f\in F({\mathcal M})$ we have:
\begin{equation*}
0 = \sum_{\overrightarrow{e}\in \overrightarrow{E}(v)} w_{\overrightarrow{e}}
 = \sum_{\overrightarrow{e}\in \overrightarrow{E}(e)} w_{\overrightarrow{e}}
 = \sum_{\overrightarrow{e}\in \overrightarrow{E}(f)} w_{\overrightarrow{e}}.
\end{equation*}
We denote by $Harm({\mathcal M})$ the vector space of harmonic vectors
of ${\mathcal M}$.
If ${\mathcal M}$ is a periodic map then we denote by
$Harm^{per}({\mathcal M})$ the vector space of periodic harmonic vector.
For any edge $e=(v_1,v_2)$ we have two directed edges
$(v_1, e)$ and $(v_2, e)$ and the relation $w_{(v_1,e)} + w_{(v_2,e)}=0$.
Thus we can define $w_e=|w_{(v_1,e)}| = |w_{(v_2,e)}|$.
When drawing plane maps and their harmonic vectors,
it is convenient to choose an orientation on each
edge and to write down the value of the harmonic vector next to it.

\begin{figure}
\begin{center}
\epsfig{height=48mm, file=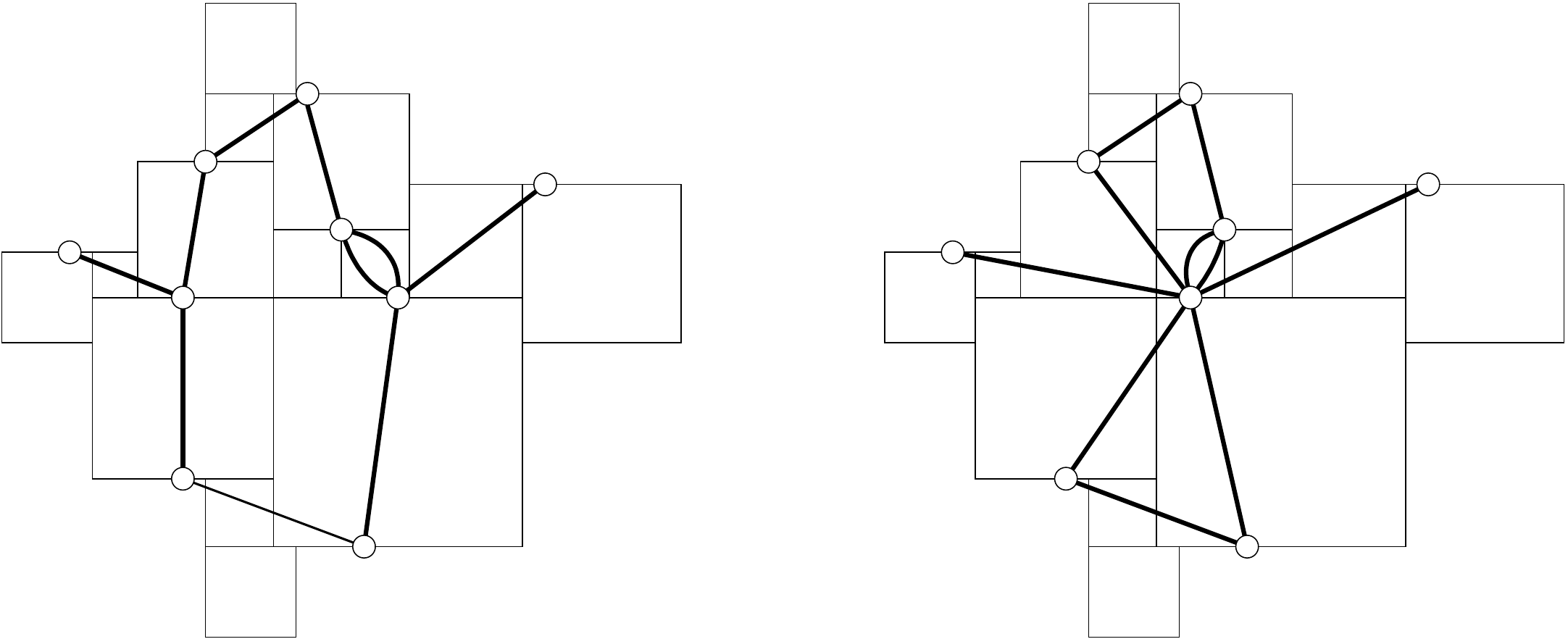}\par
\end{center}
\caption{A singular vertex and the two possible graphs ${\mathcal M}$ locally around it according to the choice of vertical and horizontal faces}
\label{ExampleSingularVertexAmbiguity}
\end{figure}

We now construct to an edge finite square tiling ${\mathcal T}$ a 
plane map ${\mathcal M}$ and a vector $w\in Harm({\mathcal M})$.
In this correspondence, every square corresponds to an edge.
Faces of ${\mathcal M}$ correspond to horizontal faces of ${\mathcal T}$
and vertices of ${\mathcal M}$ to vertical faces of ${\mathcal T}$.
For regular square tilings the map ${\mathcal M}$ is uniquely defined.
However, if ${\mathcal T}$ has a singular vertex then there are
two ways to define the horizontal and vertical faces and this affects
the map ${\mathcal M}$: either the singular vertex is contained in one
horizontal face and two vertical faces from above and below or it is
contained in one vertical face and two horizontal faces from left and
right (See in Figure \ref{ExampleSingularVertexAmbiguity} an example
of this ambiguity).
Note that all the plane tilings by squares considered in \cite{grunbaum}
are regular.

\begin{theorem}\label{FundamentalThOneWay}
For an edge finite plane square tiling ${\mathcal T}$, there exists
a plane map ${\mathcal M}$ and a harmonic vector $w\in Harm({\mathcal M})$ such that:
\begin{enumerate}
\item[(i)] Every square of ${\mathcal T}$ corresponds to an edge $e$ of ${\mathcal M}$
with $w_{e}$ being the size of the square.
\item[(ii)] For a choice of horizontal and vertical faces in
${\mathcal T}$ the map ${\mathcal M}$ is uniquely defined; its
vertices correspond to horizontal faces of ${\mathcal T}$ and
its faces to vertical faces of ${\mathcal T}$.
\item[(iii)] If ${\mathcal T}$ is periodic and the choice of horizontal
and vertical faces is done in a periodic way then the map ${\mathcal M}$ and
the vector $w$ are also periodic.
\item[(iv)] If ${\mathcal T}$ is regular then ${\mathcal M}$ is determined uniquely.
In that case, if ${\mathcal T}$ is rotated by $90$ degrees then ${\mathcal M}$ is
changed into its dual ${\mathcal M}^{*}$.
\end{enumerate}
\end{theorem}
\proof Let us take a choice of horizontal and vertical faces in ${\mathcal T}$.
We put a vertex in the middle of every horizontal face.
Also, for every square of ${\mathcal T}$ we put an edge between 
the upper horizontal face and the lower horizontal face, see
Figure \ref{ExampleSingularVertexAmbiguity} for an example.
Clearly, the vertical faces of ${\mathcal T}$ correspond to
faces of ${\mathcal M}$ and (ii) follows.
If $C$ is a square of ${\mathcal T}$ then there are vertices
$v$, $v'$ corresponding to the horizontal faces in which the
upper, respectively lower, edge of $C$ are contained.
We orient $\overrightarrow{e}$ from $v$ to $v'$ and set
$w_{\overrightarrow{e}}=d$ with $d$ the size of $C$.
For the reversed directed edge, we take the opposite value.
One checks that the vector $w$ thus defined is harmonic and (i) follows.
(iii) is clear since the above process is uniquely defined once the
horizontal and vertical faces are chosen.
If ${\mathcal T}$ is regular then there is only one possible choice of
horizontal and vertical faces. Furthermore, in that case rotating by
$90$ degrees corresponds to taking the dual and (iv) follows. \qed

We now want to define a square tiling of the plane from a plane map and
an harmonic vector.
Before that we state a lemma that is essential to this construction:
\begin{lemma}\label{CriticalSignLemma}
If ${\mathcal M}$ is a plane map and $w\in Harm({\mathcal M})$ then:
\begin{enumerate}
\item[(i)] If $v$ is a vertex of ${\mathcal M}$ and $e_1$, \dots, $e_N$ are the directed
edges coming from $v$ cyclically ordered then up to rotation the sign shape
is $(+^{a},0^{n}, -^{b}, 0^{m})$ with $a + n + b + m=N$
\item[(ii)] If $f$ is a face of ${\mathcal M}$ and $e_1$, \dots, $e_N$ are the directed
edges around $f$ cyclically ordered then up to rotation the sign shape
is $(+^{a},0^{n}, -^{b}, 0^{m})$ with $a + n + b + m=N$
\end{enumerate}

\end{lemma}
\proof This property is established in Lemma 4.11 of \cite{tutte} 
in the context of dissection of rectangles by squares.
One notices that the proof given there does
not depend on the hypothesis of conductance equal to $1$.

Now for a given plane map ${\mathcal M}$ and $w\in Harm({\mathcal M})$, let us select a vertex $v$
of ${\mathcal M}$.
We can find a block ${\mathcal B}$ of faces that contains $v$ and cut all
edges outside of ${\mathcal B}$. Now by adding edges with the right
conductance around ${\mathcal B}$ one can close the set of circuit and get
a set of wires with the current coming from a vertex $A$ and arriving at $B$.
Thus we can apply Lemma 4.11 in \cite{tutte} and conclude. \qed

We can encode dissection of rectangle by rectangles by a similar formalism
as the one of \cite{tutte}. In the correspondence, rectangles of side length
$a$, $b$ are associated to wires of conductance $a/b$.
Thus the spirit of the above lemma is that a finite set of squares can
be completed to a rectangle by adding some rectangles.

For a plane map ${\mathcal M}$ the {\em medial map} $Med({\mathcal M})$
is the map obtained by putting a vertex on every edge of ${\mathcal M}$
and having two edges adjacent if and only if they share a face and
a vertex. The resulting map $Med({\mathcal M})$ is $4$-valent and
we have $Med({\mathcal M})=Med({\mathcal M}^*)$.
The dual $(Med({\mathcal M}))^{*}$ is bipartite and its bipartite
components correspond to ${\mathcal M}$ and ${\mathcal M}^*$.

For a general map ${\mathcal M}$ and $w\in Harm({\mathcal M})$ we cannot guarantee that we will
get square tiling of the plane. Let us call an {\em unextendible square
packing} ${\mathcal T}$, a plane square packing such that around every
square there are only squares and no free space.
An example of unextendible square packing which is not a tiling is obtained
by taking the subdivision of the upper half plane $\RR_+\times \RR$
into unit squares.
Then one subdivides the squares in $[0,1]\times \RR$ into squares of
sides $1/2$ and then the squares in $[0,1/2]\times \RR$ into squares of
sides $1/4$. Continuing in this process one gets an unextendible square
packing which is not a tiling.

\begin{theorem}\label{ReconstructionResult}
(i) If ${\mathcal M}$ is a plane map and $w\in Harm({\mathcal M})$, $w\not=0$ then there exists
an unextendible square packing ${\mathcal T}({\mathcal M},w)$ such that by applying
the construction of Theorem \ref{FundamentalThOneWay} to the
corresponding choice of horizontal and vertical faces we get back
${\mathcal M}$ and $w$.

(ii) The unextendible square packing ${\mathcal T}({\mathcal M},w)$ is unique
up to translations and is simply connected.

\end{theorem}
\proof By applying Lemma \ref{CriticalSignLemma} we get that
the signs around a vertex or a face of ${\mathcal M}$ are organized in
a block of $+$ followed by a block of $-$ with possibly $0$ between them.
This gives a simple way to reconstruct the square tiling from
${\mathcal M}$ and $w$. Simply take one edge $e$ of ${\mathcal M}$ and assign it to
a square $C$ of opposite corners $(0,0)$ and $(w_e,w_e)$.
Since $e$ belongs to two vertices this allows to define the horizontal
faces over $C$ and below it.
The edge $e$ is also contained in two vertical faces and those
define vertical edges on the left and right of $C$.
Let us call two edges of ${\mathcal M}$ {\em adjacent}
if they share a vertex and a face.
So, once the position of $C$ is chosen the position of the squares
corresponding to edges adjacent to $e$ is also set.
After having put the squares around $C$ we can continue in this
way for other edges and so we get the uniqueness of ${\mathcal T}({\mathcal M},w)$
up to the position of the initial point.
But in order to prove the existence we need something else.
We need to prove that by defining the position iteratively in the
graph $Med({\mathcal M})$ we always get coherent positions.

Let us take a cycle ${\mathcal C}=\{e_1, \dots, e_N=e_1\}$ of edges with $e_i$
adjacent to $e_{i+1}$. We associate a square $C_i$ to every edge
$e_i$ according to the above defined procedure and we need to prove
that $C_1=C_{N}$.
Without loss of generality we can assume that ${\mathcal C}$ does not
self-intersect.
Let us now denote by $F_1$, \dots, $F_M$ the faces inside
of ${\mathcal C}$.
By applying the relations satisfied by $w$ around vertices,
we can modify ${\mathcal C}$
so that around any vertex of ${\mathcal C}$ the edges turn inward of it.
If we change the cycle ${\mathcal C}$ 
by removing one face $F_i$ which touches the boundary
then we do not change the position of the square.
This is because the vector $w$ is harmonic and
the sum of edge weights is zero around every face.
Thus we can remove faces $F_i$ one by one until one is reduced with the
trivial cycle $\{ e_1\}$ which proves that $C_{N}=C_1$.
Assertion (i) and (ii) follow. \qed

\section{Homology groups}\label{HomologyFormalism}

In this section we consider the space of harmonic form, first in the case
of a periodic plane graph:

\begin{theorem}\label{PeriodicStructureTh}
(i) For a periodic plane map ${\mathcal M}$ and $w\in Harm^{per}({\mathcal M})$ the unextendible
square packing ${\mathcal T}({\mathcal M}, w)$ is a periodic square tiling of the plane.

(ii) For a periodic plane graph ${\mathcal M}$, the space $Harm^{per}({\mathcal M})$ is of dimension $2$ and is isomorphic to the first homology group of a $2$-dimensional torus.
\end{theorem}
\proof Assertion (i) follows from the uniqueness condition of Theorem 
\ref{ReconstructionResult}. The fact that it is a tiling follows from the
fact that $|w(e)|\geq C$ with $C>0$ or $|w(e)|=0$. This prevents the
formation of accumulation points and guarantees that ${\mathcal T}({\mathcal M}, w)$ is a tiling.

Denote by $L=\ZZ v_1 + \ZZ v_2$ the group of translation of ${\mathcal M}$
with $v_1$ not collinear to $v_2$.
The torus $\TT={\mathcal M}/L$ is a map with
vertices, edges and faces, i.e. it is a cell complex.
Moreover, the space $Harm^{per}({\mathcal M})$ is isomorphic to $Harm(\TT)$.

For an edge $e=(v_1, v_2)$, the condition
$\sum_{\overrightarrow{e}\in \overrightarrow{E}(e)} w_{\overrightarrow{e}}=0$ is rewritten
as $w_{(v_1,e)} + w_{(v_2,e)} = 0$. This simply means that if we choose
an orientation on every edge then the space $Harm(\TT)$ is defined over
$\RR^{E(\TT)}$.
We can thus build a chain complex:
\begin{equation*}
\RR^{F(\TT)} \mathop{\longrightarrow}\limits^{d_1} \RR^{E(\TT)} \mathop{\longrightarrow}\limits^{d_2} \RR^{V(\TT)}.
\end{equation*}
If one sets the standard scalar products on $\RR^{V(\TT)}$, $\RR^{E(\TT)}$
and $\RR^{F(\TT)}$ then we can define the adjoint $d_1^{*}$ and $d_2^{*}$.
The space $Harm(\TT)$ is then $\Kert(d_2)\cap \Kert(d_1^{*})$ and
by Hodge Theory (See, for example, Theorem 7.55 in \cite{nakahara})
it is identified with the first homology group $H_1(\TT, \RR)$ with
coefficients in $\RR$.
Since $\TT$ is a torus, $H_1(\TT, \RR)$ is of dimension $2$. \qed

\begin{figure}
\begin{center}
\epsfig{height=52mm, file=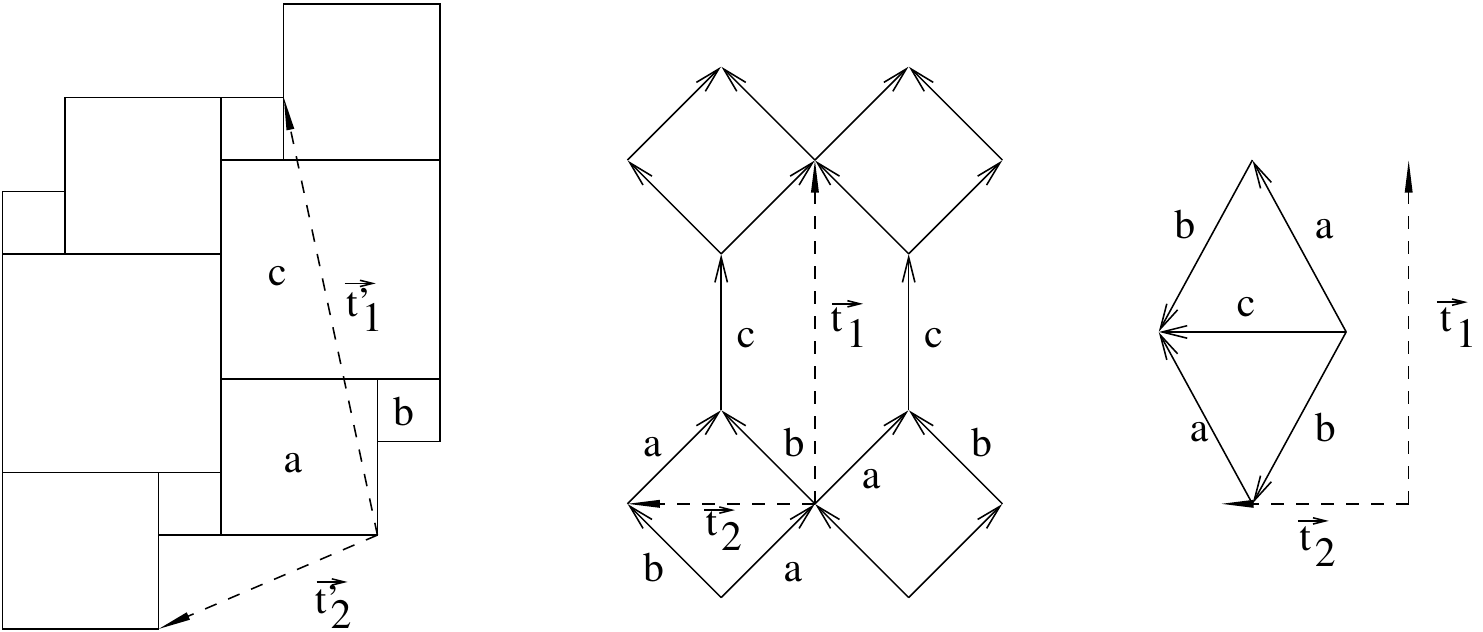}\par

\end{center}
\caption{A periodic square tiling with translation vectors $t'_1=(b-a, a+b+c)$ and $t'_2=(-a-b, a-b)$. See also the associated plane map ${\mathcal M}$ and its dual ${\mathcal M}^*$}
\label{ExampleSquareTorus}
\end{figure}

As a consequence of (ii) we know that if ${\mathcal M}$ is periodic along a lattice
$L$ and $w$ is an harmonic vector periodic along a lattice $L'\subset L$ then $w$ is periodic along $L$ as well.
If ${\mathcal M}$ is periodic then let us write the lattice $L$
of translations as $L=\ZZ \overrightarrow{t}_1 + \ZZ \overrightarrow{t}_2$.
The translations along $\overrightarrow{t}_i$ correspond to a cycle
$\gamma_i$ in the quotient ${\mathcal M} / L$.
Let us write $\gamma_i$ as a cycle of vertices
\begin{equation*}
\{v_i^1, v_i^{2},  \dots,  v_i^{l}\} \mbox{~with~} v_i^l = v_i^1 + \overrightarrow{t}_i
\end{equation*}
and the pair $e_i^k=(v_i^k, v_i^{k+1})$ being edges.
We can now define the integral of $w\in Harm^{per}({\mathcal M})$ as
\begin{equation*}
\int_{\gamma_i} w = \sum_{k=1}^{l-1} w_{(v_i^k, e_i^k)}.
\end{equation*}
The square tiling ${\mathcal T}({\mathcal M}, w)$ is periodic and we can write its lattice of translation as $L'=\ZZ \overrightarrow{t}'_1 + \ZZ \overrightarrow{t}'_2$.
The $y$ component of $\overrightarrow{t}'_i$ is $\int_{\gamma_i} w$.
For the $x$ component, one must proceed differently.
The harmonic space $Harm^{per}({\mathcal M}^*)$ is isomorphic
to $Harm^{per}({\mathcal M})$ and one can realize
the isomorphism explicitly.
It suffices to turn the oriented edges of ${\mathcal M}$ by $90$ degrees
and assign to them the same edge weight as in ${\mathcal M}$.
If $w^{*}$ is the induced harmonic vector and $\gamma_i^{*}$ are the induced
cycles then the $x$ component of $\overrightarrow{t}'_i$ is equal to
$\int_{\gamma^*_i} w^*$. See in Figure \ref{ExampleSquareTorus} a complete
example.

\begin{figure}
\begin{center}
\begin{minipage}[b]{6.2cm}
\centering
\epsfig{height=50mm, file=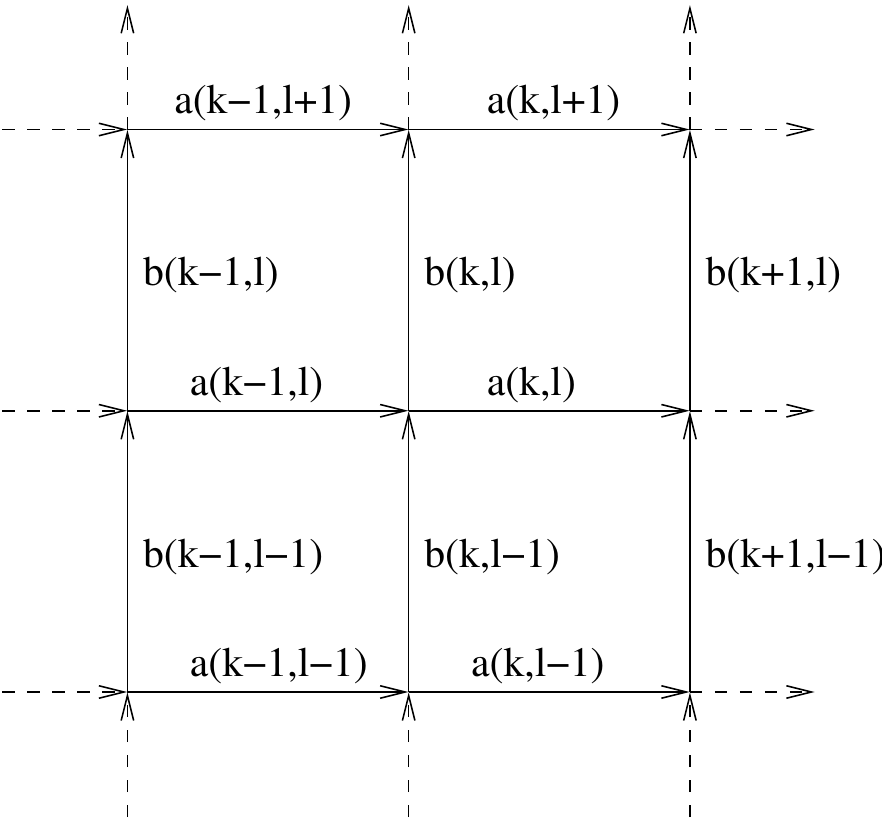}\par
\end{minipage}

\end{center}
\caption{The edge weight of the map $\ZZ^2$}
\label{EdgeWeightZ2}
\end{figure}

Theorem \ref{PeriodicStructureTh}.(ii) is classical and corresponds
to the identification of harmonic forms on a compact manifold with
real homology classes.
However, it is well known that Hodge Theory collapses in the non-compact
case. This allows us to build many interesting square tilings:

\begin{theorem}
There exists a square tiling of the plane with tiles of size all different
such that the size of the tiles grows linearly with the distance from
the origin.
\end{theorem}
\proof Let us consider the plane map $\ZZ^2$ formed by a family of squares.
Let us assign a value $a(k,l)$, $b(k,l)$ to horizontal, respectively
vertical edges according to Figure \ref{EdgeWeightZ2}.
The fact that this defines an harmonic vector corresponds
to the following equations for $k,l\in \ZZ$:
\begin{equation*}
\left\{\begin{array}{rcl}
a(k,l) - a(k-1,l) + b(k,l)   - b(k,l-1) &=& 0,\\
a(k,l) + b(k+1,l) - a(k,l+1) - b(k,l)   &=& 0.
\end{array}\right.
\end{equation*}
The two solutions below
\begin{equation*}
\begin{array}{rcl}
w_1 &=& \{a(k,l)=1  \mbox{~and~} b(k,l)=0\},\\
w_2 &=& \{a(k,l)=0  \mbox{~and~} b(k,l)=1\}
\end{array}
\end{equation*}
form a basis of $Harm^{per}(\ZZ^2)$.
But $Harm(\ZZ^2)$ is not reduced to it.
The next most interesting harmonic vectors are:
\begin{equation*}
\begin{array}{rcl}
w_3 &=& \{a(k,l)=l  \mbox{~and~} b(k,l)=k\}\\
w_4 &=& \{a(k,l)=k  \mbox{~and~} b(k,l)=-l\}
\end{array}
\end{equation*}
Now if we set the harmonic vector
\begin{equation*}
w = (1/4) w_1 + (1/3) w_2 + w_3 + \sqrt{2} w_4
\end{equation*}
then one sees that $|w(e)|=|w(e')|$ implies $e=e'$ and that the growth
of $|w(e)|$ is linear with the distance from the origin.
Since $w_e$ goes to infinity on each horizontal and vertical direction
of $\ZZ^2$, we get that the unextendible square packing ${\mathcal T}(\ZZ^2, w)$ is actually a tiling. \qed

The above theorem could be seen as a counterexample to a problem
posed by C. Pomerance and mentioned in \cite{grunbaum}:
does there exist a square tiling of the plane with all squares
of different size such that after ordering the sizes
the growth is less than exponential? However the question is not
explicitly stated there and if one imposes that square sizes
are integral then we do not have the answer to the question.

\section{Types of square tilings}\label{TypeSquareTiling}

Given a periodic plane map ${\mathcal M}$ and $w\in Harm^{per}({\mathcal M})$ we have a corresponding
periodic square plane tiling ${\mathcal T}({\mathcal M}, w)$.
By Theorem \ref{PeriodicStructureTh} we can find a basis $(w^1, w^2)$
of $Harm^{per}({\mathcal M})$.
So, $w=\alpha_1 w^1 + \alpha_2 w^2$ and this allows us to define the
corresponding $Sq$-domains:

\begin{definition}
Given a periodic plane map ${\mathcal M}$ and a basis $(w^1, w^2)$ of $Harm^{per}({\mathcal M})$
we define:

(i) $Plane({\mathcal M}, \{w^1, w^2\})$ to be the plane of the parameters $(\alpha_1, \alpha_2)$.

(ii) For an edge $e$ of ${\mathcal M}$, we choose an oriented edge $\overrightarrow{e}\in \overrightarrow{E}(e)$ and define a line $L(e)$ to be the set of $(\alpha_1, \alpha_2)\in Plane({\mathcal M}, \{w^1, w^2\})$ such that $\alpha_1 w^1_{\overrightarrow{e}} + \alpha_2 w^2_{\overrightarrow{e}}=0$.

(iii) A $Sq$-domain to be a connected set $D$ of $Plane({\mathcal M}, \{w^1, w^2\})$ which is
delimited by lines $L(e)$
such that
$\alpha_1 w^1_{\overrightarrow{e}} + \alpha_2 w^2_{\overrightarrow{e}}$ is of specified sign $\epsilon(e)$ for $\alpha=(\alpha_1, \alpha_2)\in D$.
\end{definition}

Given ${\mathcal M}$ and $(w^1, w^2)$ it is very easy to determine all possible
$Sq$-domains. For every edge $e$ of ${\mathcal M}$ we associate the line
$\alpha_1 w^1_{\overrightarrow{e}} + \alpha_2 w^2_{\overrightarrow{e}} = 0$ in the plane $(\alpha_1, \alpha_2)$.
All such lines pass by the origin and so the plane is partitioned in
a number of angular sectors which are exactly the $Sq$-domains.

Take now a point near a line $L(e)$. If one goes towards $L(e)$
then the length $w(e)$ of the corresponding square vanish.
After one passes through $L(e)$, the square reappears
but in a different position.
In Figure \ref{FamousEightUnilateralEquitransitive} the eight unilateral
plane tilings with three sizes of faces are listed. They are all regular
and we give their corresponding plane maps ${\mathcal M}$.
One notices that the square tilings (2), (3) and (4) are associated to
the same periodic plane map.

\section{Symmetry and classification considerations}\label{PossibleSymmetry}

Suppose that we have a periodic plane map ${\mathcal M}$ with
a group $\Gamma$.
Any $f\in \Gamma$ defines an action $f_{*}$ on $Harm^{per}({\mathcal M})$.
The induced action of $\Gamma$ is actually the action of the point
group of $\Gamma$.

For a square tiling, rotation axis of order $3$ or $6$ are forbidden 
and this leaves $p1$, $p2$ and $p4$ as possible rotation groups.
However, we do not know which of the $17$ wallpaper groups can occur
as group of a periodic square tiling.

Let us now see how groups of ${\mathcal M}$ and their square tiling
${\mathcal T}({\mathcal M}, w)$ are related:

\begin{proposition}
(i) If $f$ is a symmetry of a periodic plane map ${\mathcal M}$ and $w\in Harm^{per}({\mathcal M})$ the $f$ induces a symmetry of ${\mathcal T}({\mathcal M}, w)$ if and only if $f_*(w)=\pm w$.

(ii) If ${\mathcal M}$ has symmetry $p2$ then all square tilings ${\mathcal T}({\mathcal M}, w)$
also have symmetry $p2$ at least.

(iii) If $w$ is chosen at random in $Plane({\mathcal M}, \{w^1, w^2\})$
then ${\mathcal T}({\mathcal M}, w)$ has symmetry $p1$, $p2$ or $p4$.
\end{proposition}
\proof (i) is quite clear. The symmetry $p2$ is the antipodality so of course
it satisfies condition (i) and maps itself to the tiling.
For a random $w\in Harm^{per}({\mathcal M})$ no reflection can preserve $w$ and so the symmetry of a generic square tiling is $p1$, $p2$ or $p4$. \qed

In order to explain the symmetry $p4$ one has to introduce the notion of
self-duality.
If one looks at the definition of the harmonic space $Harm({\mathcal M})$ one remarks that
vertices and faces play identical roles.
A map is called {\em self-dual} if the dual ${\mathcal M}^*$ is isomorphic to ${\mathcal M}$.
For the square tiling ${\mathcal T}({\mathcal M}, w)$ taking the dual corresponds to rotating by $90$ degrees.
Thus if a regular square tiling has a $4$-fold axis then its
corresponding plane map ${\mathcal M}$ is necessarily self-dual.

\begin{figure}
\begin{center}
\begin{minipage}[b]{5.3cm}
\centering
\resizebox{50mm}{!}{\rotatebox{0}{\includegraphics[bb=125 234 536 572, clip]{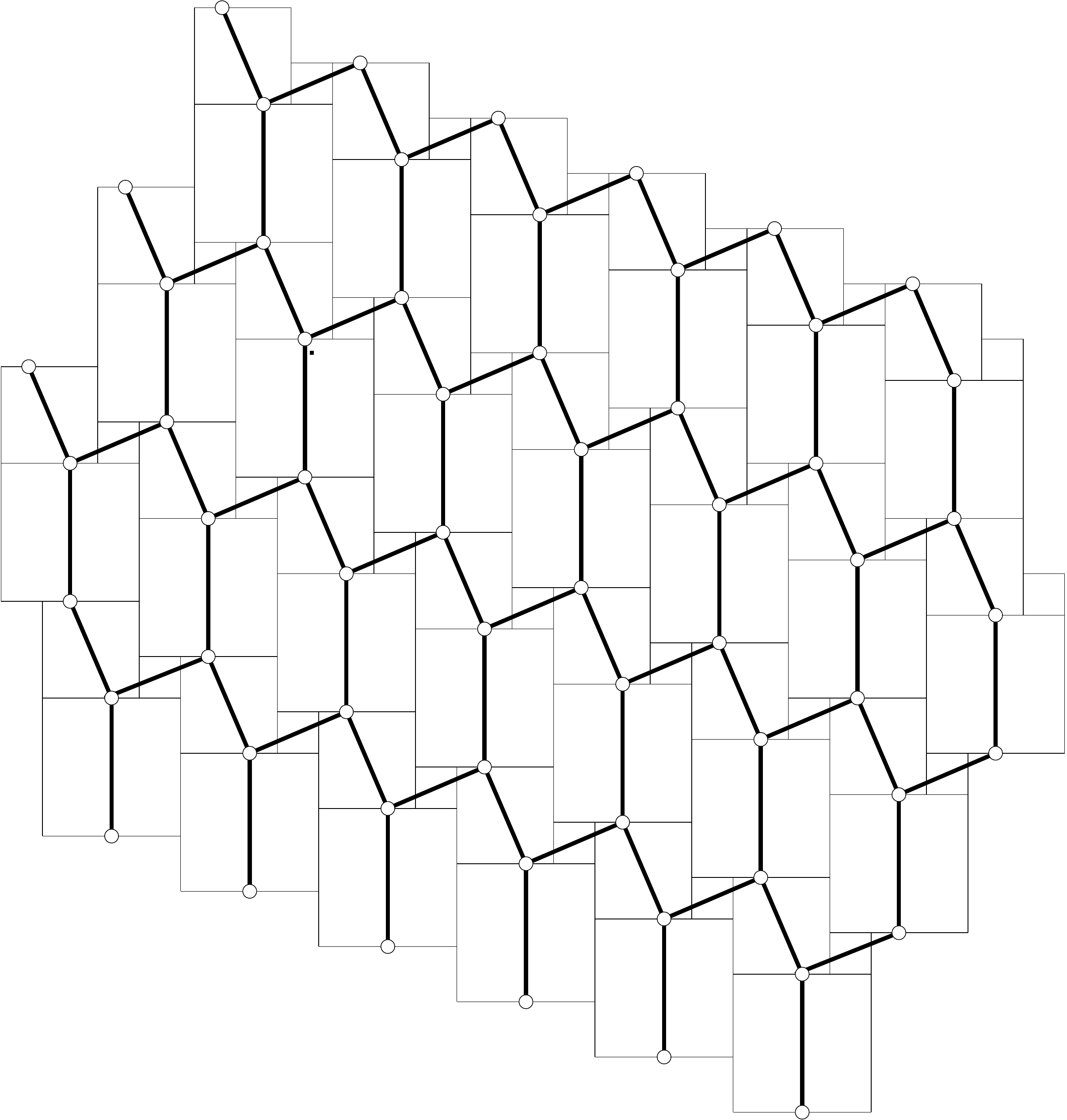}}}\par
(1): $p2$
\end{minipage}
\begin{minipage}[b]{5.3cm}
\centering
\resizebox{50mm}{!}{\rotatebox{0}{\includegraphics[bb=183 68 539 349, clip]{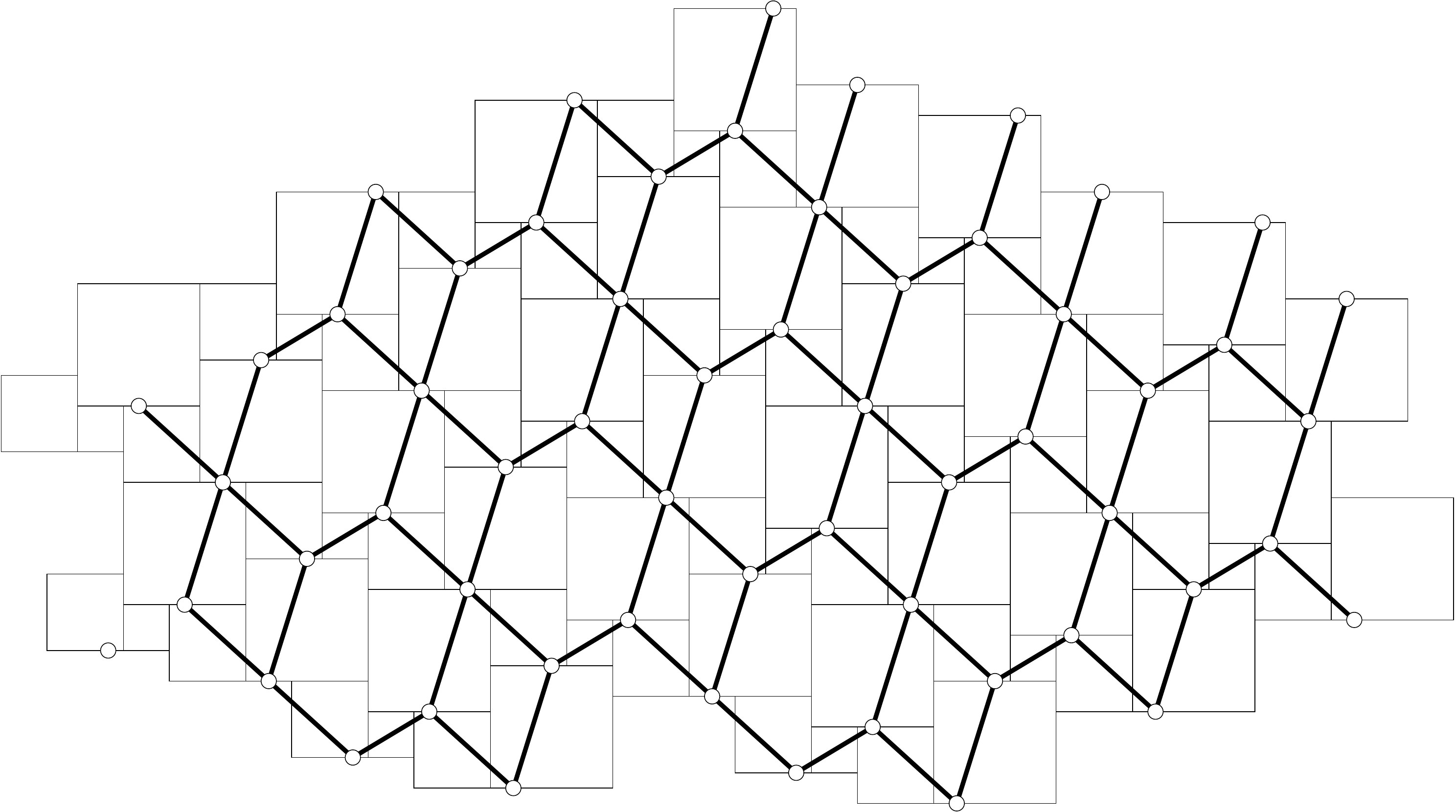}}}\par
(2): $p2$
\end{minipage}
\begin{minipage}[b]{5.3cm}
\centering
\resizebox{50mm}{!}{\rotatebox{0}{\includegraphics[bb=104 57 513 346, clip]{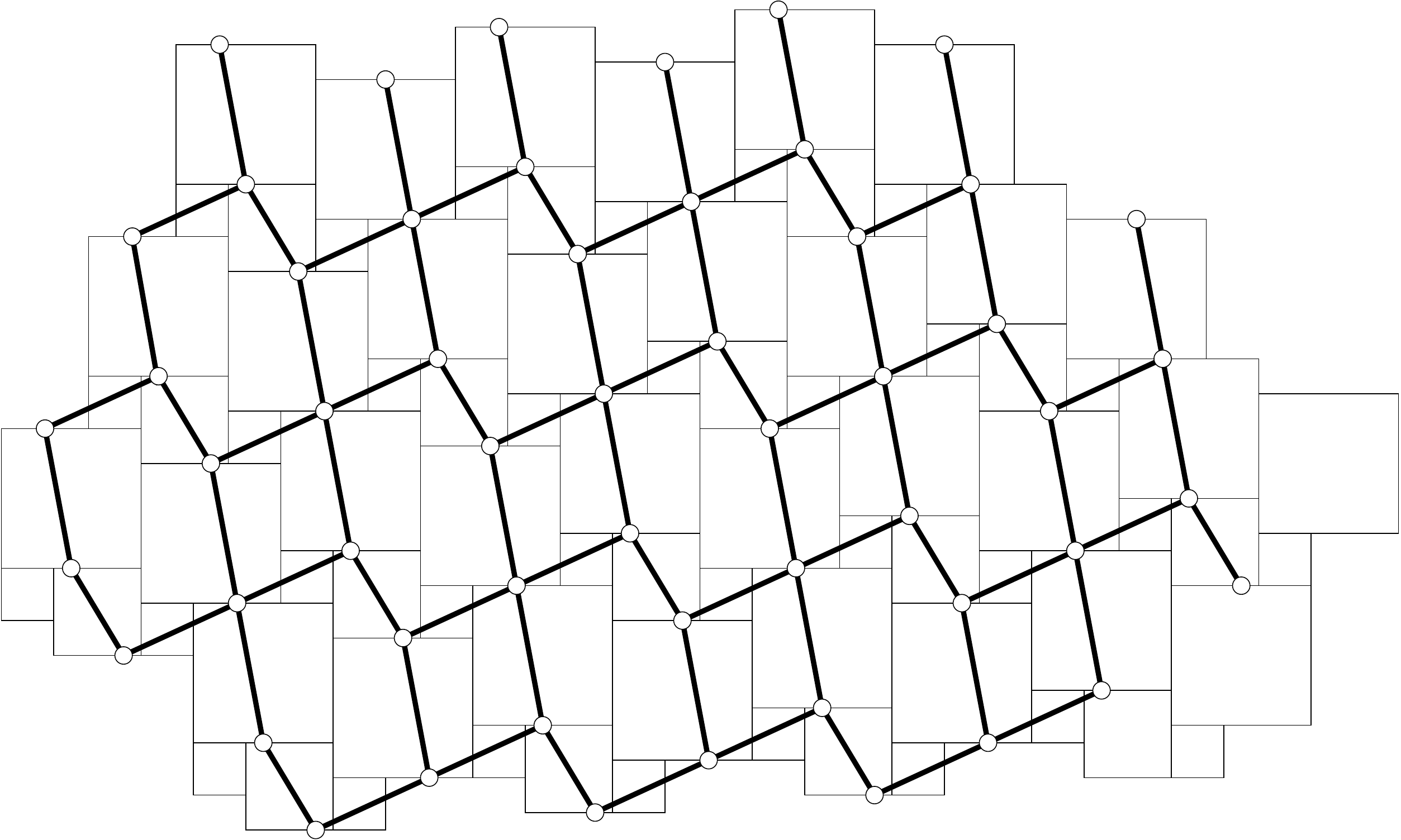}}}\par
(3): $p2$
\end{minipage}
\begin{minipage}[b]{5.3cm}
\centering
\resizebox{50mm}{!}{\rotatebox{0}{\includegraphics[bb=88 167 430 406, clip]{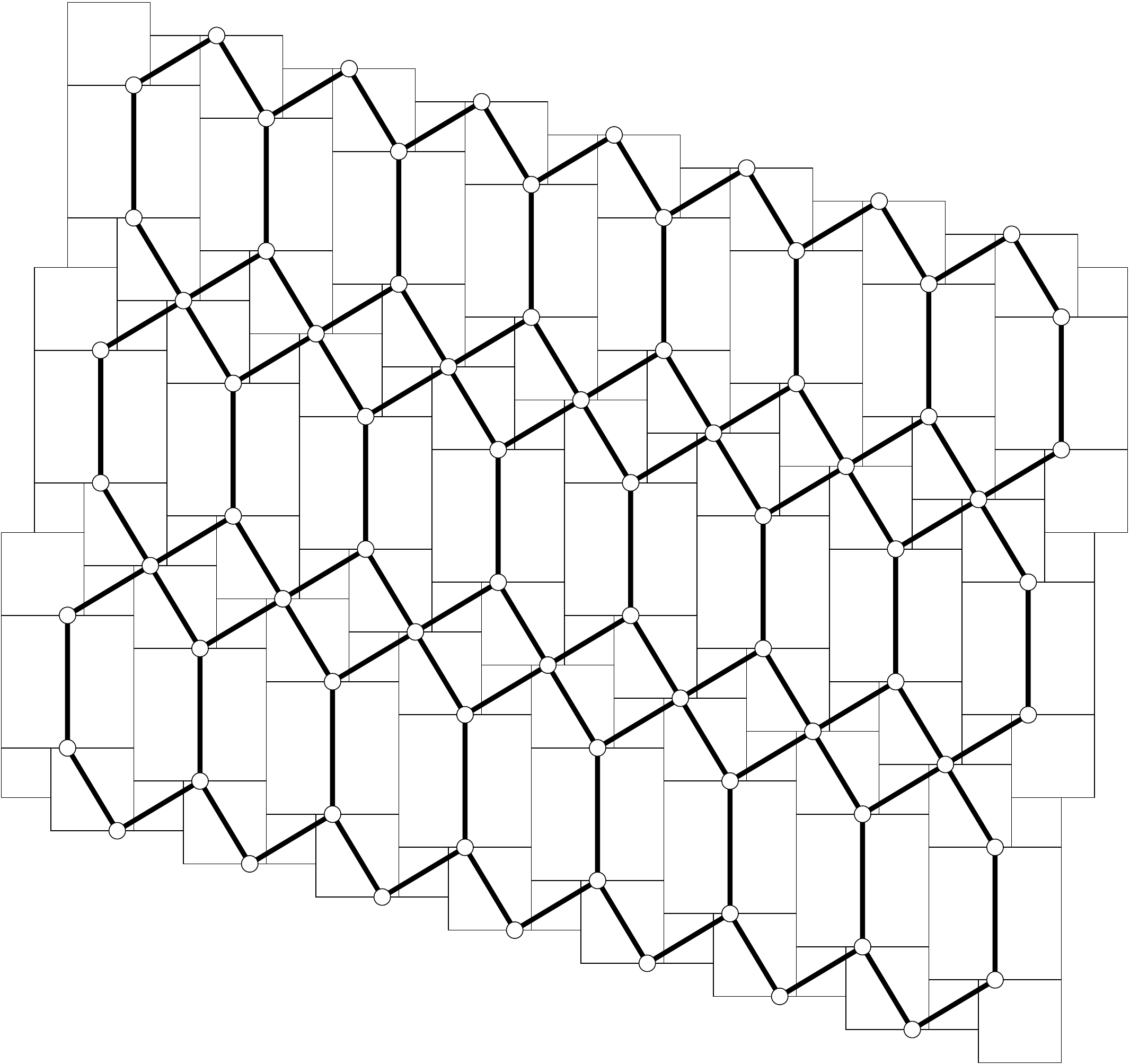}}}\par
(4): $p2$
\end{minipage}
\begin{minipage}[b]{5.3cm}
\centering
\resizebox{50mm}{!}{\rotatebox{0}{\includegraphics[bb=133 94 541 405, clip]{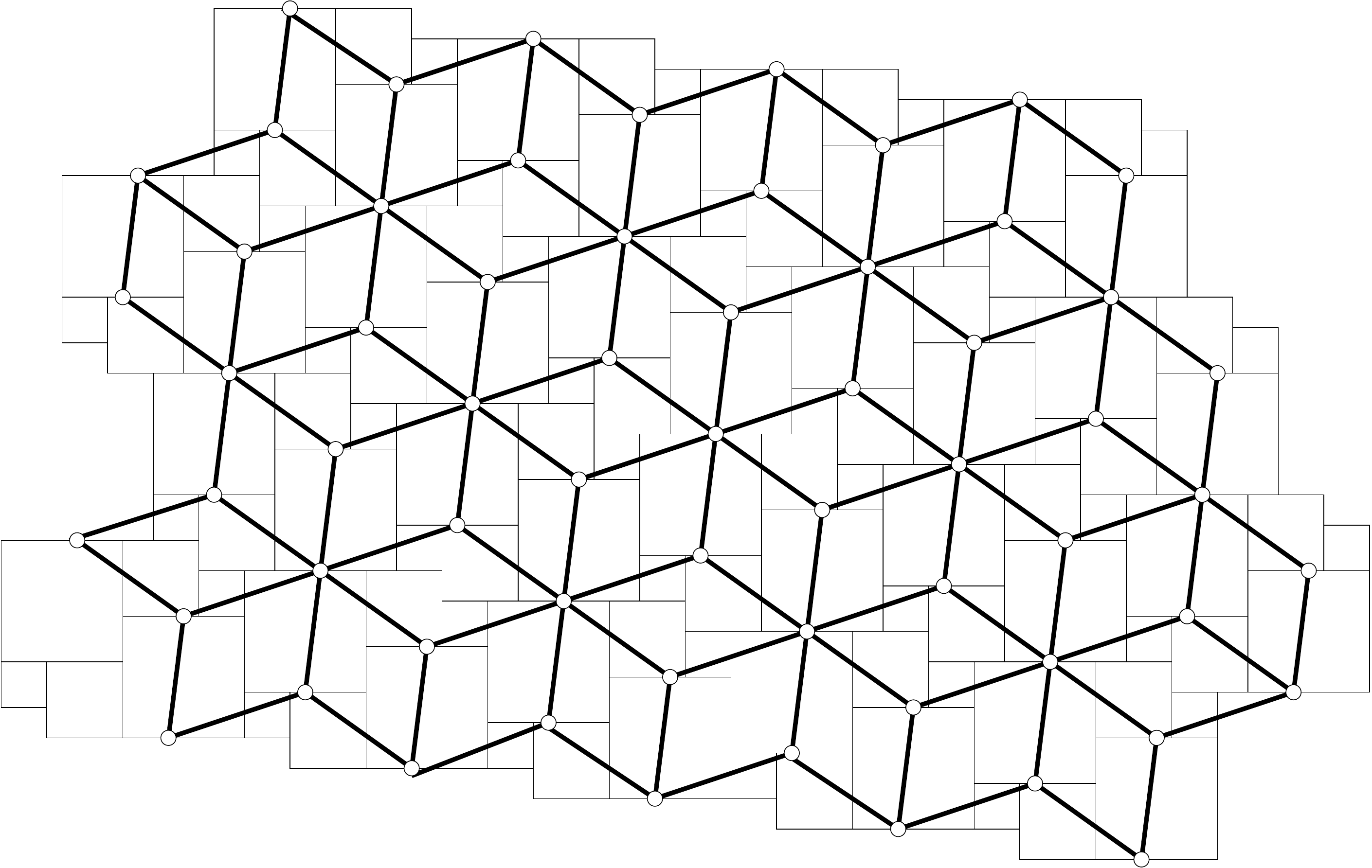}}}\par
(5): $p2$
\end{minipage}
\begin{minipage}[b]{5.3cm}
\centering
\resizebox{50mm}{!}{\rotatebox{0}{\includegraphics[bb=88 72 629 483, clip]{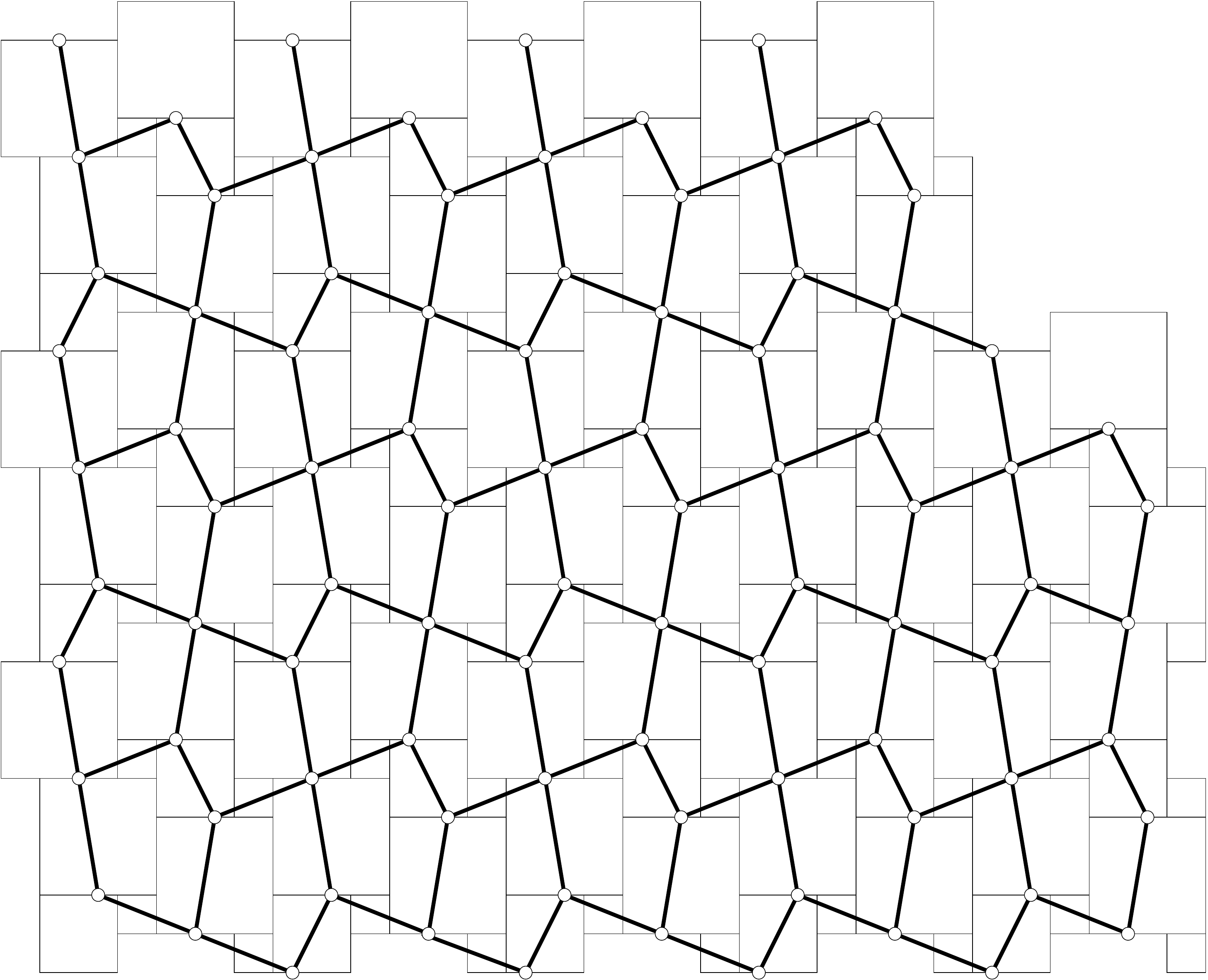}}}\par
(6): $pgg$
\end{minipage}
\begin{minipage}[b]{5.3cm}
\centering
\resizebox{50mm}{!}{\rotatebox{0}{\includegraphics[bb=295 259 735 643, clip]{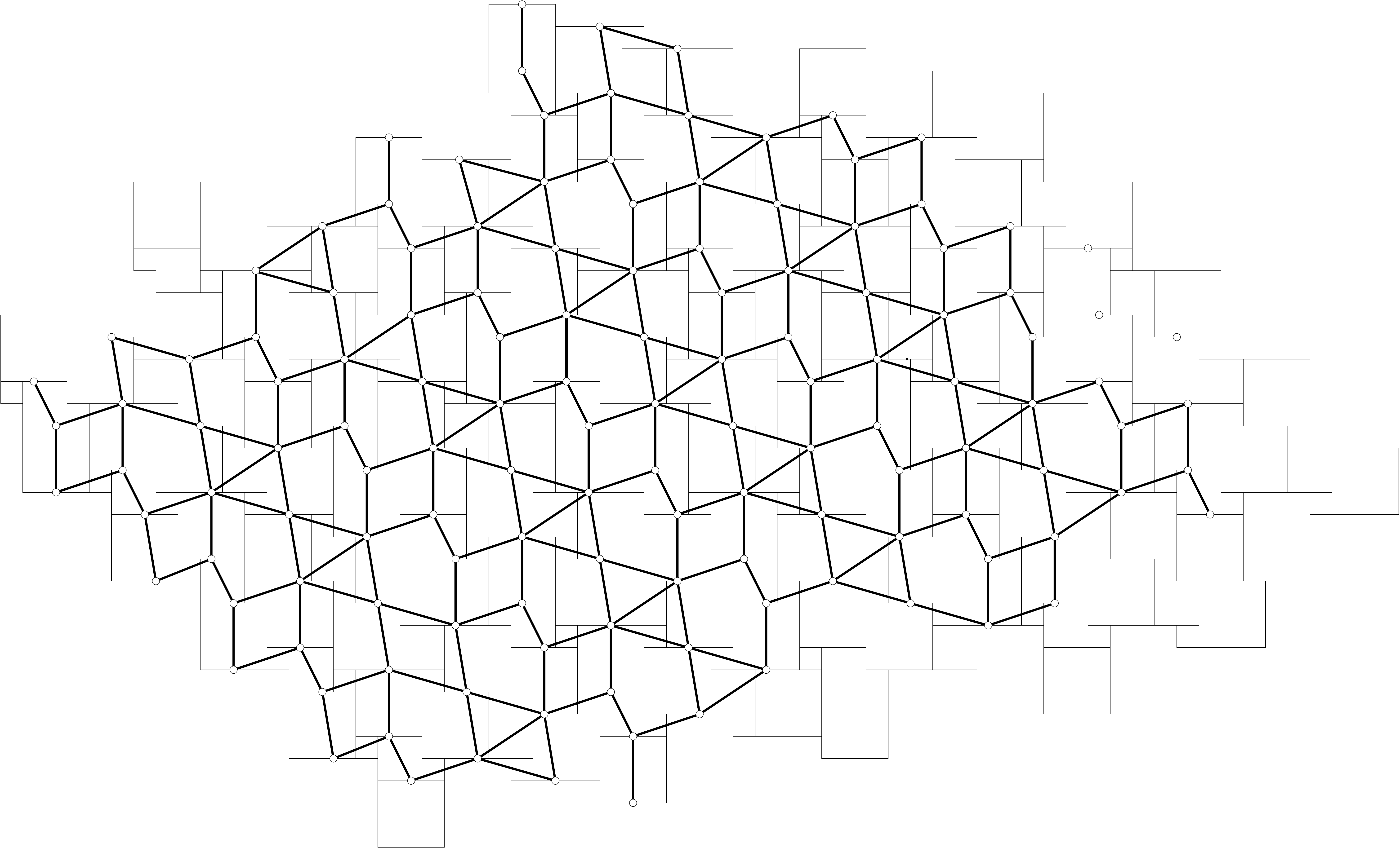}}}\par
(7): $pgg$
\end{minipage}
\begin{minipage}[b]{5.3cm}
\centering
\resizebox{50mm}{!}{\rotatebox{0}{\includegraphics[bb=264 308 740 722, clip]{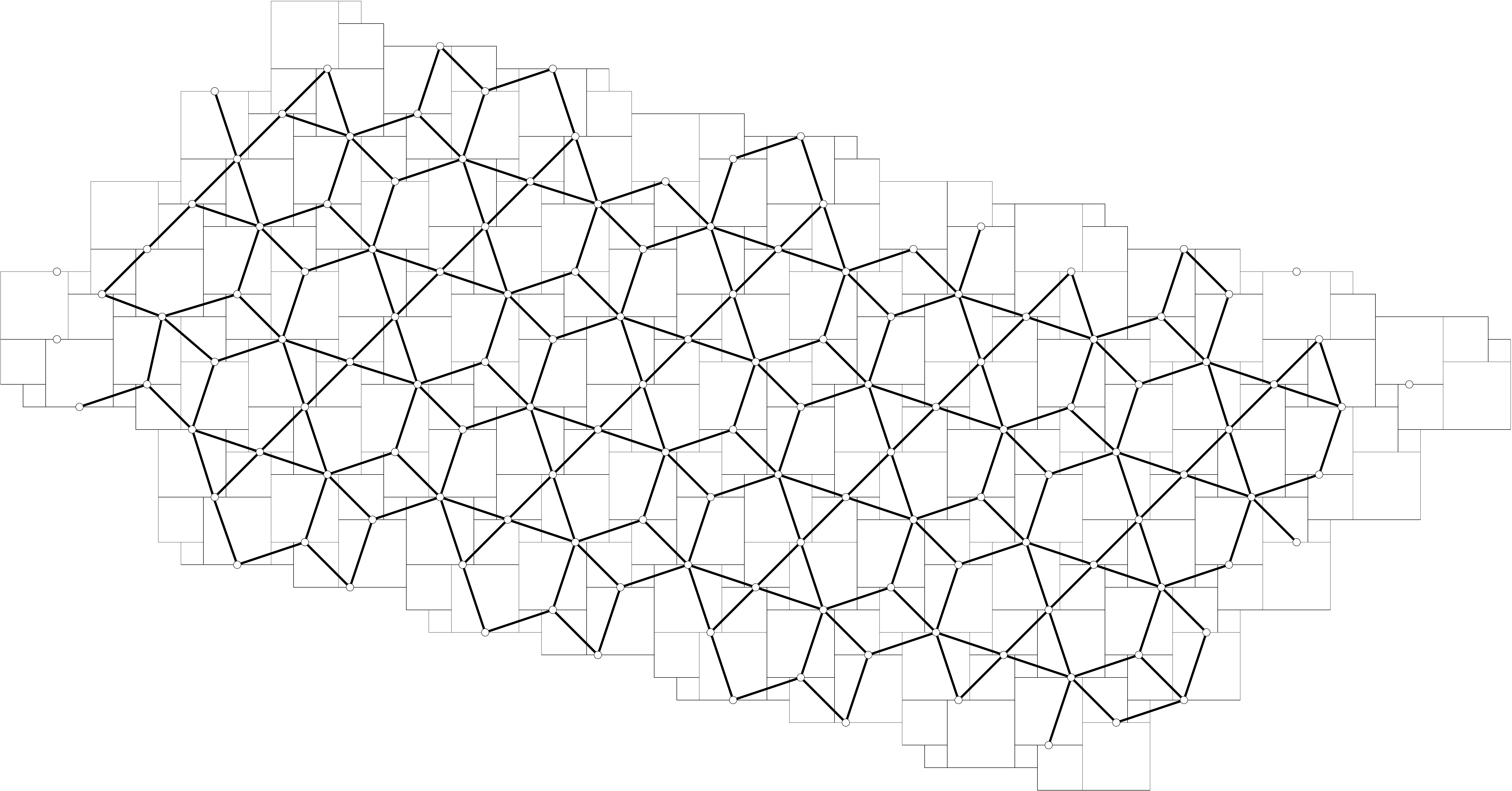}}}\par
(8): $pgg$
\end{minipage}

\end{center}
\caption{The $8$ square plane tiling with $3$ orbits of squares and their corresponding plane maps}
\label{FamousEightUnilateralEquitransitive}
\end{figure}

Figure \ref{FamousEightUnilateralEquitransitive} gives for each of the eight
plane tiling with $3$ orbits of squares of different size the corresponding
plane map ${\mathcal M}$.
The first five such tilings have symmetry $p2$ and are defined by two
parameters.
The last three also have a line of symmetry and this comes from a symmetry
of the plane map and harmonic vector.
If a symmetry preserves a vector then this fixes some parameter and thus
those tilings are defined by just $1$ parameter.

We cannot use the construction that we have introduced for the enumeration
of plane square tilings with $r$ orbits of squares.
The first reason is the restriction of edge finiteness.
The second reason is that if there is a singular vertex then there is an
element of choice in the construction of the map ${\mathcal M}$, which 
greatly complexifies the enumeration.

However, if one restricts to regular square tilings, then we can design
a possible enumeration strategy. If the square tiling ${\mathcal T}$ does not
have a $4$-fold axis of symmetry but has $r$ orbits of squares 
then the corresponding map ${\mathcal M}$ has $r$ orbits of edges.
But if ${\mathcal T}$ has a $4$-fold axis of symmetry then it comes
from the self-duality of ${\mathcal M}$.
The medial operation satisfies $Med({\mathcal M})=Med({\mathcal M}^*)$
and any self-duality of ${\mathcal M}$
corresponds to an ordinary symmetry of $Med({\mathcal M})$.
Moreover, the maps ${\mathcal M}$ and ${\mathcal M}^*$ can be
reconstructed easily from $Med({\mathcal M})$ and this was already used
for enumeration of self-dual maps in \cite{selfdual}.

Thus we have to enumerate the toroidal maps having $r$ orbits of vertices
in order to get the regular square tilings with $r$ orbits of squares
by applying the technique of Delaney symbol \cite{2isohedral} and the methods
of exhaustive enumeration \cite{brinkmann}.
This should actually be possible for small values of $r$ like $r=4$, $5$.
However, as noted in \cite{grunbaum} what would be really nice is if
the existing square plane tilings of
Figure \ref{FamousEightUnilateralEquitransitive} were actually used in city
squares.

\end{document}